\documentclass{article}
\usepackage{amsmath,amssymb}  

\newcommand{\reals}{\mathbb{R}} 
\newcommand{\complex}{\mathbb{C}}

\begin{document}

\makeatletter
\def\det{\mathop{\operator@font det}\nolimits}
\def\adj{\mathop{\operator@font adj}\nolimits}
\def\rank{\mathop{\operator@font rank}\nolimits}
\makeatother

\def\dtau{\mathrm{d}\tau}                                                    
\def\exp{\mathrm{e}}  


\title{Generalization of the Matrix Determinant Lemma and its application to the controllability of single input control systems\thanks{%
Mathematics Subject Classifications: 15A24, 93B05, 93C05.}}
\date{{\small  August 30, 2017}}
\author{Robert Vrabel\thanks{%
Slovak University of Technology in Bratislava, Faculty of Materials Science and Technology in Trnava,  J. Bottu 25,  917 01 Trnava,  Slovakia (robert.vrabel@stuba.sk)}}
\maketitle

\begin{abstract}
Linear control theory provides a rich source of inspiration and motivation for development in the matrix theory. Accordingly, in this paper, a generalization of Matrix Determinant Lemma to the finite sum of outer products of column vectors is derived and an alternative proof of one of the fundamental results in modern control theory of the linear time--invariant systems $\dot x=Ax+Bu,$ $y=Cx$ is given, namely that the state controllability is unaffected by state feedback, and even more specifically, that for the controllability matrices $\mathcal{C}$ of the single input open and closed loops the equality $\det\left(\mathcal{C}_{(A,B,C)}\right)$ $=\det\left(\mathcal{C}_{(A-BK,B,C)}\right)$ holds.
\end{abstract}

\section[Introduction and main result]{Introduction and main result}

Stability, controllability and observability are the important structural properties of dynamical systems and represent three major concepts of modern control system theory. With "stability", "controllability" and "observability", one can classify the control systems without first finding the solution in an explicit form. The last two concepts were introduced by R. Kalman in the early 1960s, see \cite{Kal1}, \cite{Kal2}, \cite{Kal3}, \cite{Kal4}. The concept stability has a longer history, first mentioned in 1892 by A. M. Lyapunov in his Doctoral Dissertation \cite{Lyapunov}. As we will see later in the paper, these three concepts are closely related. 

\centerline{}

Let us consider a linear time--invariant (LTI) control system modeled by the state equation 
\begin{equation}
\dot x(t)=Ax(t)+Bu(t), \quad x(0)=x_0 \label{1}
\end{equation}
describing the time evolution of the $n$-dimensional state $x(t)$ of the system and the algebraic output equation 
\begin{equation}
y(t)=Cx(t). \label{2}
\end{equation} 
The control $u(t)$ is a $m$-dimensional vector function of time which must be chosen to make the system 
behave in a desired manner. $A,$ $B$ and $C$ are $n\times n,$  $n\times m$ and $q\times n$ constant matrices, respectively, and $x_0$ is an initial state of the system. 

\centerline{}

In order to ensure clarity and consistency of presentation, let us define these three key concepts.  

\centerline{}

{\bf Stability.} The LTI system, described by state equation (\ref{1}) is said to be
(globally) asymptotically stable if the homogeneous response (that is, $u(t)\equiv 0$) of the state vector $x(t)$ returns to the origin $x=0$ of the state space from any initial condition $x_0$ as time $t\longrightarrow\infty,$ or in terms of the state transition matrix
\[
\lim\limits_{t\to\infty}\exp^{At}x(0)=0
\]
for any $x(0).$

Equivalently, for any perturbations of initial condition $x(0)=x_0$ the perturbed solutions are attracted by the solution of (\ref{1}) for $t\longrightarrow\infty.$   

\centerline{}

{\bf Controllability.} The LTI system (\ref{1}), (\ref{2}) or the triple $(A,B,C)$ is said to be state controllable if for any initial state $x(0)=x_0$ and any final state $x_{_T},$ there exists a control input $\bar u(t)$ that transfers $x_0$ to $x_{_T}$ in a finite time $T.$ In other words
\[
x_{_T}=\exp^{AT}x_0+\int\limits_0^T\exp^{A(T-\tau)}B\bar u(\tau)\dtau.
\]
Otherwise $(A,B,C)$ is said to be uncontrollable.

\centerline{}

{\bf Observability.} The LTI system, is said to be observable if any initial state $x(0)$ can be uniquely determined from the knowledge of the output $y(t)$ and the input $u(t)$ on the interval $[0,t_1]$ for some $t_1>0.$ Otherwise $(A,B,C)$ is said to be unobservable. 

\centerline{}

The simple algebraic conditions can be given for the asymptotic stability, controllability and observability of the LTI control system (\ref{1}), (\ref{2}):

\smallskip

THEOREM~1. A necessary and sufficient condition for $(A,B,C)$ to be asymptotically stable is that every eigenvalue of $A$ has a strictly negative real part.

\smallskip

THEOREM~2. (Kalman) A necessary and sufficient condition for $(A,B,C)$ to be controllable is
\[
\rank\mathcal{C}_{(A,B,C)}=:\rank\left[ B\ \vdots\ AB\ \vdots\ A^2B\ \vdots\ \cdots\ \vdots\ A^{n-1}B\right]=n.
\] 

\smallskip

THEOREM~3. (Kalman) A necessary and sufficient condition for $(A,B,C)$ to be observable is
\[
\rank\mathcal{O}_{(A,B,C)}=:\rank\left[ C^T\ \vdots\ A^TC^T\ \vdots\ (A^T)^2C^T\ \vdots\ \cdots\ \vdots\ (A^T)^{n-1}C^T\right]^T=n.
\] 
A superscript $T$ denotes the matrix transpose operation.

\centerline{}

Consider the state feedback control $u(t)=-Kx(t) + r(t),$ where $K$ is an $m \times n$ gain matrix and $r(t)$ is an $m-$dimensional external input. The closed loop system dynamics is given by
\[
\dot x(t)=(A - BK)x(t)+Br(t).
\]
A fundamental result of linear control theory is that the five following conditions are equivalent: 
\begin{itemize}
\item[(i)] the triple $(A,B,C)$ is controllable;
\item[(ii)] the triple $\left(A^T,C^T,B^T\right)$ is observable \ (controllability--observability duality);
\item[(iii)] $\rank\mathcal{C}_{(A,B,C)}=n$ \ (Kalman test);
\item[(iv)] $\rank(A-\lambda I_n \ \vdots\ B)=n,$ for every $\lambda\in\complex$ \ (Belovich-Popov-Hautus test);
\item[(v)] for every $\alpha_0,\alpha_1,\dots,\alpha_{n-1}\in\reals$ there exists a matrix $K\in\reals^{m\times n}$ such that
\[
\chi_{_{A-BK}}(\lambda)=\lambda^n+\alpha_{n-1}\lambda^{n-1}+\dots+\alpha_1\lambda+\alpha_0, 
\]
where $\chi_{_{A-BK}}$ is the characteristic polynomial of the matrix $A-BK$ \  (pole placement problem or state feedback stabilization problem).
\end{itemize}

The equivalence $(i)\Longleftrightarrow (ii)$ can often be used to go from the results on controllability to ones on observability, and vice versa.  

\centerline{}

Another important property is the fact that controllability is unaffected by state feedback.

\smallskip

THEOREM~4. $(A,B,C)$ is controllable if and only if $(A-BK,B,C)$ is controllable for all $K$ of dimension $m\times n.$

\smallskip

PROOF. The statement follows from the Belovich-Popov-Hautus test and the matrix identity
\[
\left[(A-BK)-\lambda I_n \ \vdots\ B\right]=\left[A-\lambda I_n \ \vdots\ B\right]
{\setlength\arraycolsep{2pt}
\left[\begin{array}{rr}
  I_n\ & 0   \\
  -K \ & I_n
\end{array} \right], 
}
\]
or using another matrix identity \cite[p.~181]{BanksTran}
\begin{equation}
\mathcal{C}_{(A-BK,B,C)}=\mathcal{C}_{(A,B,C)}\mathcal{T}, \label{3} 
\end{equation}
where
\[
\mathcal{C}_{(A-BK,B,C)}=\left[ B\ \vdots\ (A-BK)B\ \vdots\ (A-BK)^2B\ \vdots\ \cdots\ \vdots\ (A-BK)^{n-1}B\right],
\]
\[
\mathcal{C}_{(A,B,C)}=\left[ B\ \vdots\ AB\ \vdots\ A^2B\ \vdots\ \cdots\ \vdots\ A^{n-1}B\right]
\]
and
\small{
\[
\mathcal{T}={\setlength\arraycolsep{2pt}
\left[
\begin{array}{cccccc}
  I_m    \ & -KB     \ & -K(A-BK)B\ & -K(A-BK)^2B\ & \cdots\ & -K(A-BK)^{n-2}B   \\
  0      \ & I_m     \ & -KB      \ & -K(A-BK)B  \ & \cdots\ & -K(A-BK)^{n-3}B   \\
  0      \ & 0       \ & I_m      \ & -KB        \ & \cdots\ & -K(A-BK)^{n-4}B   \\
  0      \ & 0       \ & 0        \ & I_m        \ & \cdots\ & -K(A-BK)^{n-5}B   \\
  \vdots \ & \vdots  \ & \vdots   \ & \vdots     \ & \ddots \ & \vdots           \\
  0      \ & 0       \ & 0        \ & 0          \ & \cdots\ & I_m
\end{array} 
\right].
}
\]
}

\smallskip

In this paper we prove the stronger statement of THEOREM~4 for the case when the controllability matrices $\mathcal{C}_{(A,B,C)}$ and $\mathcal{C}_{(A-BK,B,C)}$ are the square matrices, namely: 

\smallskip

THEOREM~5. For $m=1$ and $n\geq2$
\[
\det\left(\mathcal{C}_{(A,B,C)}\right)=\det\left(\mathcal{C}_{(A-BK,B,C)}\right).
\]

\smallskip

EXAMPLE~1. As a simple but illustrative example, consider the single input LTI control system
\[
{\setlength\arraycolsep{2pt}
\left[\begin{array}{c}
  \dot x_1   \\
  \dot x_2
\end{array} \right]=
\left[\begin{array}{rr}
  1\ & 1   \\
  1\ & 2
\end{array} \right]\left[\begin{array}{c}
  x_1   \\
  x_2
\end{array} \right]+\left[\begin{array}{c}
  1   \\
  0
\end{array} \right][u],\qquad y=Cx.
}
\] 
The controllability matrix of open-loop system is
\[
\mathcal{C}_{(A,B,C)}=
{\setlength\arraycolsep{2pt}
\left[\begin{array}{rr}
  1\ & 1   \\
  0\ & 1
\end{array} \right],
}\ \mathrm{that\ is,} \ \det\left(\mathcal{C}_{(A,B,C)}\right)=1.
\]
Now for an arbitrarily chosen state feedback gain matrix $K=\left(k_1\ k_2\right)$ we get for closed-loop system that
\[
\mathcal{C}_{(A-BK,B,C)}=
{\setlength\arraycolsep{2pt}
\left[\begin{array}{rc}
  1\ & 1-k_1   \\
  0\ & 1
\end{array} \right],
\ \mathrm{that\ is,}\ \det\left(\mathcal{C}_{(A-BK,B,C)}\right)=1.
}
\]

\smallskip

REMARK 1. Notice that the statement of THEOREM~5 immediately follows from the identity (3), bearing in the mind that $\det\left(\mathcal{T}\right)=1.$ The aim of this paper is to provide an alternative proof of this theorem based on the Generalized Matrix Determinant Lemma to gain deeper insight into structure of the controllability matrices of open and closed loops.

\section{Proof of the main result} 

We precede the proof of THEOREM~5 by following new lemma which is a generalization of the well-known Matrix Determinant Lemma (MDL) representing an important analytical tool in the matrix theory, theory of optimal control, etc.  For the wider context of MDL, see e.~g. \cite{Harville}. Its application in the theory of control of biological systems can be found in \cite{CaMidHu} and the problem of multidimensional root finding by using MDL is studied in \cite{NoferiniTownsend}. By employing the MDL, the optimal sensor placement problem to achieve optimal measurements was investigated in \cite{SongPasionLhommeOldenburg}. 

\smallskip

LEMMA~1~{\bf (Generalized Matrix Determinant Lemma)}. 
Suppose $H$ is a square matrix of dimension $n$ and $u_i, v_i$ are the $n\times1$ column vectors, $i=1,\dots,k.$ Then for every $k\geq1$ we have the equality
\begin{equation}
\det\left(H+\Delta_k\right)=\det(H)+\sum\limits_{i=1}^k v_i^T\adj(H+\Delta_{i-1})u_i, \label{4}
\end{equation}
where
\[
\Delta_i = \left\{ \begin{array}{ll}
n\times n\ \mathrm{zero\ matrix} & \mathrm{for}\ i=0, \\
\sum\limits_{j=1}^i u_jv_j^T & \mathrm{for}\ i=1,\dots,k.
\end{array} \right.
\]

\smallskip

REMARK~2. For $k=1$ and an invertible matrix $H$ we obtain the classical MDL for the outer product of two vectors $u_1$ and $v_1.$

\smallskip

PROOF~OF~LEMMA~1. We use the induction principle to prove LEMMA~1. We define as a predicate $P(k)$ the statement of LEMMA~1.

{\bf Step 1:} We prove that the formula (4) is true for $k=1,$ that is,
\[
P(k=1): \quad \det\left(H+u_1v_1^T \right)=\det(H)+v_1^T\adj(H)u_1.
\]
First, let us assume that a matrix $H$ is invertible. From the matrix identity
\[
{\setlength\arraycolsep{2pt}
\left[
\begin{array}{cr}
  I_n\ & 0   \\
  v_1^T \ & 1
\end{array} 
\right]
\left[
\begin{array}{cr}
  I_n+u_1v_1^T\ & u_1   \\
  0 \ & 1
\end{array} 
\right]
\left[
\begin{array}{cr}
  I_n\ & 0   \\
  -v_1^T \ & 1
\end{array} 
\right]=
\left[
\begin{array}{rc}
  I_n\ & u_1   \\
  0 \ & 1+v_1^Tu_1
\end{array} 
\right]
}
\]
we obtain that
\[
\det\left(I_n+u_1v_1^T\right)=1+v_1^Tu_1.
\]
Hence
{\setlength\arraycolsep{1pt}
\begin{eqnarray*}
\det\left(H+u_1v_1^T\right)&=&\det(H)\det\left(I_n+(H^{-1}u_1)v_1^T\right)=\det(H)\left(1+v_1^T(H^{-1}u_1)\right) \nonumber \\
&=&\det(H)+v_1^T\adj(H)u_1.
\end{eqnarray*}}
Now let $\det(H)=0.$ Let us consider a small perturbation of $H$ in the form $H+\epsilon I_n.$ The $\det(H+\epsilon I_n)$ is a polynomial in $\epsilon$ which has at most $n$ roots on the real axis. Thus there exists $\epsilon_0$ such that the matrices $H+\epsilon I_n$ are the invertible matrices for all $\epsilon\in(0,\epsilon_0)$ and so
\[
\det\left((H+\epsilon I_n)+u_1v_1^T \right)=\det(H+\epsilon I_n)+v_1^T\adj(H+\epsilon I_n)u_1.
\]
Now in the limit for $\epsilon\rightarrow 0^+,$ taking into consideration that the polynomials on the both sides of the last equality are continuous functions, we obtain the statement $P(1).$ This completes the proof of $P(1)$ for an arbitrary square matrix $H.$ 

{\bf Step 2:} {\it (Proof that an implication $P(k=s)\Longrightarrow P(k=s+1)$ is true).} The induction hypothesis is 
that (4) is true for some $k=s\geq1.$ We have
{\setlength\arraycolsep{1pt}
\begin{eqnarray*}
\det\left(H+\Delta_{s+1}\right)&=&\det\left(\left[H+\Delta_{s}\right]+u_{s+1}v_{s+1}^T \right) \nonumber \\
&=&\det\left(H+\Delta_{s}\right)+v_{s+1}^T\adj\left(H+\Delta_{s}\right)u_{s+1} \nonumber \\
&=&\det(H)+\sum\limits_{i=1}^s v_i^T\adj(H+\Delta_{i-1})u_i+v_{s+1}^T\adj\left(H+\Delta_{s}\right)u_{s+1} \nonumber \\
&=&\det(H)+\sum\limits_{i=1}^{s+1} v_i^T\adj(H+\Delta_{i-1})u_i.
\end{eqnarray*}}
Thus (4) is true for all $k\geq1.$ 

\smallskip

REMARK~3.  From the just proved lemma it follow some useful corollaries:
\begin{itemize} 
\item[(a)] The product $UV^T$ of the matrices $U=\left[u_1\, \vdots\, u_2\, \vdots \cdots \vdots\, u_r\right]$ and $V=\left[ v_1\, \vdots\, v_2\,\vdots \cdots \vdots\, v_r\right],$ where $u_i$ and $v_i$ are $n\times 1$ column vectors, $i=1,\dots,r$ ($r\geq1$) may be expressed in the form of the sum of outer products, $UV^T=\sum\limits_{n=1}^r u_iv_i^T.$ Thus from (4) we have the matrix determinant identity
\[
\det\left(H+UV^T\right)=\det(H)+\sum\limits_{i=1}^r v_i^T\adj(H+\Delta_{i-1})u_i,
\]
where $H$ is an arbitrary $n\times n$ matrix and $\Delta_i$ are defined in LEMMA~1. This equality can be used also for deriving some results if the matrix $U$ or/and $V$ have a special form (see the proof of THEOREM~5 below);
\item[(b)] For $H=0$ (the zero matrix) we obtain a generalized formula for the determinant of the product of two (in general non-square) matrices 
\[
\det\left(UV^T\right)=\sum\limits_{i=1}^r v_i^T\adj(\Delta_{i-1})u_i.
\]
For $r=1$ we get the obvious fact that the matrix formed by the outer product of two vectors has a determinant equal to zero.
\end{itemize}

\smallskip

PROOF OF THEOREM~5. Let us denote the column vectors $(A-BK)^iB-A^iB$ by $\beta_i,$ $i=1,\dots,n-1.$ Now we rewrite $\mathcal{C}_{(A-BK,B,C)}$ as the following sum:
{\setlength\arraycolsep{1pt}
\begin{eqnarray*}
\mathcal{C}_{(A-BK,B,C)}&=&\left[ B\ \vdots\ (A-BK)B\ \vdots\ (A-BK)^2B\ \vdots\ \cdots\ \vdots\ (A-BK)^{n-1}B\right] \nonumber \\
&=&\mathcal{C}_{(A,B,C)}+\left[ 0\ \vdots\ \beta_1\ \vdots\ \cdots\ \vdots\ \beta_{n-1}\right]=\mathcal{C}_{(A,B,C)}+ \underbrace{\left[\beta_1\right]}_{u_1}\underbrace{[0\ 1\ 0\ \dots\ 0]}_{v_1^T} \nonumber \\
&+&\underbrace{\left[\beta_2\right]}_{u_2}\underbrace{[0\ 0\ 1\ \dots\ 0]}_{v_2^T}+\dots + \underbrace{\left[\beta_{n-1}\right]}_{u_{n-1}}\underbrace{[0 \ 0\ 0\ \dots\ 1]}_{v_{n-1}^T}\nonumber \\
&=&\mathcal{C}_{(A,B,C)}+\sum\limits_{i=1}^{n-1} u_iv_i^T.
\end{eqnarray*}}
Then on the basis of LEMMA~1 for $H=\mathcal{C}_{(A,B,C)}$ and $k=n-1$ we obtain the equality
\[
\det\left(\mathcal{C}_{(A-BK,B,C)}\right)=\det\left(\mathcal{C}_{(A,B,C)}\right)+\sum\limits_{i=1}^{n-1} v_i^T\adj(\mathcal{C}_{(A,B,C)}+\Delta_{i-1}) u_i.
\]
In the iterative way we show that $v_i^T\adj(\mathcal{C}_{(A,B,C)}+\Delta_{i-1}) u_i=0$ for all $i=1,\dots,n-1:$
\begin{itemize}
\item[\ \ $i=1:$] $v_1^T\adj(\mathcal{C}_{(A,B,C)}) u_1=\det(H_1),$ where $H_1$ is the matrix obtained from $\mathcal{C}_{(A,B,C)}$ by replacing the second column of $\mathcal{C}_{(A,B,C)}$ by the column $\beta_1;$
\item[$i=2:$] $v_2^T\adj(\mathcal{C}_{(A,B,C)}+u_1v_1^T) u_2=\det(H_2),$ where $H_2$ is the matrix obtained from $H_1$ by replacing the third column of $H_1$ by the column $\beta_2;$
\item[$\vdots$]
\item[$i=n-1:$] $v_{n-1}^T\adj(\mathcal{C}_{(A,B,C)}+\Delta_{n-2}) u_{n-1}=\det(H_{n-1}),$ where $H_{n-1}$ is the matrix obtained from $H_{n-2}$ by replacing the $n-$th column of $H_{n-2}$ by the column $\beta_{n-1}.$
\end{itemize}
Because  
\[
[B]\left(=\mathrm{1^{st}\ column\ in\ }H_i\right)(-KB)=-BKB=\left[\beta_1\right]\left(=\mathrm{2^{nd}\ column\ in\ }H_i\right),
\]
the determinant of $H_i$ is $0,$ for all $i=1,\dots,n-1.$
This completes the proof of THEOREM~5.

\section{Theoretical conclusions for control theory}

Let us consider that the initial state of the LTI control system (\ref{1}) is at origin. Due to Cayley-Hamilton theorem, the final states $x_{_T}$ at the time $t=T$ can be written as a linear combination of the columns of Kalman controllability matrix $\mathcal{C}_{(A,B,C)}.$ If we think of $\mathcal{C}_{(A,B,C)}$ and $\mathcal{C}_{(A-BK,B,C)}$ as representing the linear transformations, then the nonzero determinant (more precisely, the absolute value of the determinant) represents how much the linear transformation is stretching or compressing  the bounded regions in $\reals^n$ \cite{Hannah}. As follows from THEOREM~5, the "$n-$dimensional volume distortion" ratio for open and closed loop is the same. This fact may be interpreted that the basis for achieving the state from the bounded region $S_f$ of the state space $\reals^n$ at the finite time $T$ have 
equal volumes,
\[
\mathrm{Volume\ of}\ \mathcal{C}^{-1}_{(A,B,C)}(S_f)=\mathrm{Volume\ of}\ \mathcal{C}^{-1}_{(A-BK,B,C)}(S_f),
\]
for the open loops and the loops with feedback. 

\section*{Information on the copyright holder:}
The parts of this manuscript have been published as separate papers in the International Journal of Pure and Applied Mathematics, Vol. 111 (4), pp.~643-646 (2016) under the title "A Note on the Matrix Determinant Lemma" (doi: 10.12732/ijpam.v111i4.11) and  Vol. 116 (1), pp.~147-151 (2017) under the title  "Generalized Matrix Determinant Lemma and the Controllability of Single Input Control Systems" (doi: 10.12732/ijpam.v116i1.15). 

\smallskip

\section*{Acknowledgement}
This publication is the result of implementation of the project "University Scientific Park:  Campus MTF STU - CAMBO" (26220220179) supported by the Research and Development Operational Program (ERDF).

\end{document}